# Additive isotone regression

Enno Mammen[1,*] and Kyusang Yu[1,*]

*Universität Mannheim*

**This paper is dedicated to Piet Groeneboom
on the occasion of his 65th birthday**

**Abstract:** This paper is about optimal estimation of the additive components of a nonparametric, additive isotone regression model. It is shown that asymptotically up to first order, each additive component can be estimated as well as it could be by a least squares estimator if the other components were known. The algorithm for the calculation of the estimator uses backfitting. Convergence of the algorithm is shown. Finite sample properties are also compared through simulation experiments.

## 1. Introduction

In this paper we discuss nonparametric additive monotone regression models. We discuss a backfitting estimator that is based on iterative application of the pool adjacent violator algorithm to the additive components of the model. Our main result states the following oracle property. Asymptotically up to first order, each additive component is estimated as well as it would be (by a least squares estimator) if the other components were known. This goes beyond the classical finding that the estimator achieves the same rate of convergence, independently of the number of additive components. The result states that the asymptotic distribution of the estimator does not depend on the number of components.

We have two motivations for considering this model. First of all we think that this is a useful model for some applications. For a discussion of isotonic additive regression from a more applied point, see also Bacchetti [1], Morton-Jones et al. [32] and De Boer, Besten and Ter Braak [7]. But our main motivation comes from statistical theory. We think that the study of nonparametric models with several nonparametric components is not fully understood. The oracle property that is stated in this paper for additive isotone models has been shown for smoothing estimators in some other nonparametric models. This property is expected to hold if the estimation of the different nonparametric components is based on local smoothing where the localization takes place in different scales. An example are additive models of smooth functions where each localization takes place with respect to another covariate. In Mammen, Linton and Nielsen [28] the oracle property has been verified for the local linear smooth backfitting estimator. As local linear estimators, also the isotonic least squares is a local smoother. The estimator is a local average of the response variable but in contrast to local linear estimators the local neighborhood is chosen

---

*Research of this paper was supported by the Deutsche Forschungsgemeinschaft project MA 1026/7-3 in the framework of priority program SPP-1114.

[1]Department of Economics, Universität Mannheim, L 7, 3–5, 68131 Mannheim, Germany, e-mail: emammen@rumms.uni-mannheim.de; yukyusan@rumms.uni-mannheim.de

by the data. This data adaptive choice is automatically done by the least squares minimization. This understanding of isotonic least squares as a local smoother was our basic motivation to conjecture that for isotonic least squares the oracle property should hold as for local linear smooth backfitting.

It may be conjectured that the oracle property holds for a much larger class of models. In Horowitz, Klemela and Mammen [19] a general approach was introduced for applying one-dimensional nonparametric smoothers to an additive model. The procedure consists of two steps. In the first step, a fit to the additive model is constructed by using the projection approach of Mammen, Linton and Nielsen [28]. This preliminary estimator uses an undersmoothing bandwidth, so its bias terms are of asymptotically negligible higher order. In a second step, a one-dimensional smoother operates on the fitted values of the preliminary estimator. For the resulting estimator the oracle property was shown: This two step estimator is asymptotically equivalent to the estimator obtained by applying the one-dimensional smoother to a nonparametric regression model that only contains one component. It was conjectured that this result also holds in more general models where several nonparametric components enter into the model. Basically, a proof could be based on this two step procedures. The conjecture has been verified in Horowitz and Mammen [20, 22] for generalized additive models with known and with unknown link function.

The study of the oracle property goes beyond the classical analysis of rates of convergence. Rates of convergence of nonparametric estimators depend on the entropy of the nonparametric function class. If several nonparametric functions enter into the model the entropy is the sum of the entropies of the classes of the components. This implies that the resulting rate coincides with the rate of a model that only contains one nonparametric component. Thus, rate optimality can be shown for a large class of models with several nonparametric components by use of empirical process theory, see e.g. van de Geer [39]. Rate optimality for additive models was first shown in Stone [38]. This property was the basic motivation for using additive models. In contrast to a full dimensional model it allows estimation with the same rate of convergence as a one-dimensional model and avoids for this reason the curse of dimensionality. On the other hand it is a very flexible model that covers many features of the data nonparametrically. For a general class of nonparametric models with several components rate optimality is shown in Horowitz and Mammen [21].

The estimator of this paper is based on backfitting. There is now a good understanding of backfitting methods for additive models. For a detailed discussion of the basic statistical ideas see Hastie and Tibshirani [18]. The basic asymptotic theory is given in Opsomer and Ruppert [34] and Opsomer [35] for the classical backfitting and in Mammen, Linton and Nielsen [28] for the smooth backfitting. Bandwidth choice and practical implementations are discussed in Mammen and Park [29, 30] and Nielsen and Sperlich [33]. The basic difference between smooth backfitting and backfitting lies in the fact that smooth backfitting is based on a smoothed least squares criterion whereas in the classical backfitting smoothing takes place only for the updated component. The full smoothing of the smooth backfitting algorithm stabilizes the numerical and the statistical performance of the estimator. In particular this is the case for degenerated designs and for the case of many covariates as was shown in simulations by Nielsen and Sperlich [33]. In this paper we use backfitting without any smoothing. For this reason isotone additive least squares will have similar problems as classical backfitting and these problems will be even more severe because no smoothing is used at all. Smooth backfitting methods for generalized additive models were introduced in Yu, Park and Mammen [42]. Haag [15] discusses smooth backfitting for nonparametric additive diffusion models. Tests based on smooth backfitting have been considered in Haag [16] and Mammen and Sperlich [31]. Backfitting tests have been proposed in Fan and Jiang [10]. Additive regression is an example of a nonparametric model where the nonparametric function is given as a solution of an integral equation. This has been outlined in Linton and Mammen [24] and Carrasco, Florens and Renault [6] where also other examples of statistical integral equations are given. Examples are additive models where the additive components are linked as in Linton and Mammen [25] and regression models with dependent errors where an optimal transformation leads to an additive model, see Linton and Mammen [26]. The representation of estimation in additive models as solving an empirical integral equation can also be used to understand why the oracle property holds.

In this paper we verify the oracle property for additive models of isotone functions. It is shown that each additive component can be estimated with the same asymptotic accuracy as if the other components would be known. We compare the performance of a least squares backfitting estimator with a least squares isotone estimator in the oracle model where only one additive component is unknown. The backfitting estimator is based on iterative applications of isotone least squares to each additive component. Our main theoretical result is that the differences between these two estimators are of second order. This result will be given in the next section. The numerical performance of the isotone backfitting algorithm and its numerical convergence will be discussed in Section 3. Simulations for the comparison of the isotone backfitting estimator with the oracle estimator are presented in Section 4. The proofs are deferred to the Appendix.

## 2. Asymptotics for additive isotone regression

We suppose that we have i.i.d. random vectors $(Y^1, X_1^1, \ldots, X_d^1), \ldots, (Y^n, X_1^n, \ldots, X_d^n)$ and we consider the regression model

(1) $$E(Y^i | X_1^i, \ldots, X_d^i) = c + m_1(X_1^i) + \cdots + m_d(X_d^i)$$

where $m_j(\cdot)$'s are monotone functions. Without loss of generality we suppose that all functions are monotone increasing. We also assume that the covariables take values in a compact interval, $[0, 1]$, say. For identifiability we add the normalizing condition

(2) $$\int_0^1 m_j(x_j) \, dx_j = 0.$$

The least squares estimator for the regression model (1) is given as minimizer of

(3) $$\sum_{i=1}^n (Y^i - c - \mu_1(X_1^i) - \cdots - \mu_d(X_d^i))^2$$

with respect to monotone increasing functions $\mu_1, \ldots, \mu_d$ and a constant $c$ that fulfill $\int_0^1 \mu_j(x_j) \, dx_j = 0$. The resulting estimators are denoted as $\widehat{m}_1, \ldots, \widehat{m}_d$ and $\widehat{c}$.

We will compare the estimators $\widehat{m}_j$ with oracle estimators $\widehat{m}_j^{OR}$ that make use of the knowledge of $m_l$ for $l \neq j$. The oracle estimator $\widehat{m}_j^{OR}$ is given as minimizer

of

$$\sum_{i=1}^{n}(Y^i - c - \mu_j(X_1^i) - \sum_{l \neq j} m_l(X_l^i))^2$$
$$= \sum_{i=1}^{n}(m_j(X_j^i) + \varepsilon^i - c - \mu_j(X_1^i))^2$$

with respect to a monotone increasing function $\mu_j$ and a constant $c$ that fulfill $\int_0^1 \mu_j(x_j)dx_j = 0$. The resulting estimators are denoted as $\widehat{m}_j^{OR}$ and $\widehat{c}^{OR}$.

In the case $d = 1$, this gives the isotonic least squares estimator proposed by Brunk [4] which is given by

$$\widehat{m}_1(X_1^{(i)}) = \max_{s \leq i} \min_{t \geq i} \sum_{j=s}^{t} Y^{(j)}/(t - s + 1) \tag{4}$$

where $X_1^{(1)}, \ldots, X_1^{(n)}$ are the order statistics of $X^1, \ldots, X^n$ and $Y^{(j)}$ is the observation at the observed point $X_1^{(j)}$. Properties and simple computing algorithms are discussed e.g. in Barlow et al. [2] and Robertson, Wright, and Dykstra [36]. A fast way to calculate the estimator is to use the Pool Adjacent Violator Algorithm (PAVA). In the next section we discuss a backfitting algorithm for $d > 1$ that is based on iterative use of PAVA.

We now state a result for the asymptotic performance of $\widehat{m}_j$. We use the following assumptions. To have an economic notation, in the assumptions and in the proofs we denote different constants by the same symbol $C$.

(A1) The functions $m_1, \ldots, m_d$ are differentiable and their derivatives are bounded on $[0, 1]$. The functions are strictly monotone, in particular for $G(\delta) = \inf_{|u-v| \geq \delta, 1 \leq j \leq d} |m_j(v) - m_j(u)|$ it holds $G(\delta) \geq C\delta^\gamma$ for constants $C, \gamma > 0$ for all $\delta > 0$.

(A2) The $d$-dimensional vector $X^i = (X_1^i, \ldots, X_d^i)$ has compact support $[0, 1]^d$. The density $p$ of $X^i$ is bounded away from zero and infinity on $[0, 1]^d$ and it is continuous. The tuples $(X^i, Y^i)$ are i.i.d. For $j, k = 1, \ldots, d$ the density $p_{X_k, X_j}$ of $(X_k^i, X_j^i)$ fulfills the following Lipschitz condition for constants $C, \rho > 0$

$$\sup_{0 \leq u_j, u_k, v_k \leq 1} |p_{X_k, X_j}(u_k, u_j) - p_{X_k, X_j}(v_k, u_j)| \leq C|u_k - v_k|^\rho.$$

(A3) Given $X^i$ the error variables $\varepsilon^i = Y^i - c - m_1(X_1^i) - \cdots - m_d(X_d^i)$ have conditional zero mean and subexponential tails, i.e. for some $\gamma > 0$ and $C' > 0$, it holds that

$$E\left[\exp(\gamma|\varepsilon^i|)\Big|X^i\right] < C' \qquad \text{a.s.}$$

The conditional variance of $\varepsilon^i$ given $X^i = x$ is denoted by $\sigma^2(x)$. The conditional variance of $\varepsilon^i$ given $X_1^i = u_1$ is denoted by $\sigma_1^2(u_1)$. We assume that $\sigma_1^2$ is continuous at $x_1$.

These are weak smoothness conditions. We need (A3) to apply results from empirical process theory. Condition (A1) excludes the case that a function $m_j$ has flat parts. This is done for the following reason. Isotonic least squares regression produces piecewise constant estimators where for every piece the estimator is equal to the sample average of the piece. If the function is strictly monotone the pieces shrink to 0, a.s. If the function has flat parts these averages do not localize at the flat parts. But in our proof we make essential use of a localization argument. We conjecture that our oracle result that is stated below also holds for the case that there are flat parts. But we do not pursue to check this here. It is also of minor interest because at flat parts the monotone least squares estimator is of order $\widehat{m}_j^{OR} - m_j = o_P(n^{-1/3})$. Thus the oracle result $\widehat{m}_j - \widehat{m}_j^{OR} = o_P(n^{-1/3})$ then only implies that $\widehat{m}_j - m_j = o_P(n^{-1/3})$. In particular, it does not imply that $\widehat{m}_j$ and $\widehat{m}_j^{OR}$ have the same asymptotic distribution limit.

For $d = 1$ the asymptotics for $\widehat{m}_1$ are well known. Note that the estimator $\widehat{m}_1$ for $d = 1$ coincides with the oracle estimator $\widehat{m}_1^{OR}$ for $d > 1$ that is based on isotonizing $Y^i - c - m_2(X_2^i) - \cdots - m_d(X_d^i) = m_1(X_1^i) + \varepsilon^i$ in the order of the values of $X_1^i$ ($i = 1, \ldots, n$). For the oracle model (or for the model (1) with $d = 1$) the following asymptotic result holds under (A1)–(A3):

*For all $x_1 \in (0,1)$ it holds that*

$$\widehat{m}_1^{OR}(x_1) - m_1(x_1) = O_P(n^{-1/3}).$$

*Furthermore at points $x_1 \in (0,1)$ with $m_1'(x_1) > 0$, the normalized estimator*

$$n^{1/3} \frac{[2p_1(x_1)]^{1/3}}{\sigma_1(x_1)^{2/3} m_1'(x_1)^{1/3}} [\widehat{m}_1^{OR}(x_1) - m_1(x_1)]$$

*converges in distribution to the slope of the greatest convex minorant of $W(t) + t^2$, where $W$ is a two-sided Brownian motion. Here, $p_1$ is the density of $X_1^i$.*

The greatest convex minorant of a function $f$ is defined as the greatest convex function $g$ with $g \leq f$, pointwise. This result can be found, e.g. in Wright [41] and Leurgans [23]. Compare also Mammen [27]. For further results on the asymptotic law of $\widehat{m}_1^{OR}(x_1) - m_1(x_1)$, compare also Groeneboom [12, 13].

We now state our main result about the asymptotic equivalence of $\widehat{m}_j$ and $\widehat{m}_j^{OR}$.

**Theorem 1.** *Make the assumptions (A1)–(A3). Then it holds for c large enough that*

$$\sup_{n^{-1/3} \leq x_j \leq 1 - n^{-1/3}} |\widehat{m}_j(x_j) - \widehat{m}_j^{OR}(x_j)| = o_P(n^{-1/3}),$$

$$\sup_{0 \leq x_j \leq 1} |\widehat{m}_j(x_j) - \widehat{m}_j^{OR}(x_j)| = o_P(n^{-2/9} (\log n)^c)$$

The proof of Theorem 1 can be found in the Appendix. Theorem 1 and the above mentioned result on $\widehat{m}_1^{OR}$ immediately implies the following corollary.

**Corollary 1.** *Make the assumptions (A1)–(A3). For $x_1 \in (0,1)$ with $m_1'(x_1) > 0$ it holds that*

$$n^{1/3} \frac{[2p_1(x_1)]^{1/3}}{\sigma_1(x_1)^{2/3} m_1'(x_1)^{1/3}} [\widehat{m}_1(x_1) - m_1(x_1)]$$

*converges in distribution to the slope of the greatest convex minorant of $W(t) + t^2$, where $W$ is a two-sided Brownian motion.*

## 3. Algorithms for additive isotone regression

The one-dimensional isotonic least squares estimator can be regarded as a projection of the observed vector $(Y^{(1)}, \ldots, Y^{(n)})$ onto the convex cone of isotonic vectors in $\mathbb{R}^n$ with respect to the scalar product $\langle f, g \rangle \equiv \sum_{i=1}^{n} f^{(i)} g^{(i)}$ where $f \equiv (f^{(1)}, \ldots, f^{(n)})$ and $g \equiv (g^{(1)}, \ldots, g^{(n)}) \in \mathbb{R}^n$. Equivalently, we can regard it as a projection of a right continuous simple function with values $(Y^{(1)}, \ldots, Y^{(n)})$ onto the convex cone of right continuous simple monotone functions which can have jumps only at $X_1^{(i)}$'s. The projection is with respect to the $L_2$ norm defined by the empirical measure, $P_n(y, x_1)$ which gives mass $1/n$ at each observations $(Y^i, X_1^i)$. Other monotone functions $m$ with $m(X^{(i)}) = g^{(i)}$ would also solve the least square minimization.

Now, we consider the optimization problem (3). Without loss of generality, we drop the constant. Let $H_j$, $j = 1, \ldots, d$ be the sets of isotonic vectors of length $n$ or right continuous monotone simple functions which have jumps only at $X_j^i$'s with respect to the ordered $X_j$'s. It is well known that these sets are convex cones. Then, our optimization problem can be written as follows:

$$(5) \qquad \min_{g \in H_1 + \cdots + H_d} \sum_{i=1}^{n} (Y^i - g^i)^2.$$

We can rewrite (5) as an optimization problem over a product sets $H_1 \times \cdots \times H_d$. Say $\mathbf{g} = (g_1, \ldots, g_d) \in H_1 \times \cdots \times H_d$ where $g_j \in H_j$ for $j = 1, \ldots, d$. Then the minimization problem (5) can be represented as minimizing a function over a cartesian product of sets, i.e.,

$$(6) \qquad \min_{\mathbf{g} \in H_1 \times \cdots \times H_d} F(\mathbf{Y}, \mathbf{g}).$$

Here, $F(\mathbf{Y}, \mathbf{g}) = \sum_{i=1}^{n} (Y^i - g_1^i - \cdots - g_d^i)^2$.

A classical way to solve an optimization problem over product sets is a cyclic iterative procedure where at each step we minimize $F$ with respect to one $g_j \in H_j$ while keeping the other $g_k \in H_k$, $j \neq k$ fixed. That is to generate sequences $g_j^{[r]}$, $r = 1, 2, \ldots$, $j = 1, \ldots, d$, recursively such that $g_j^{[r]}$ minimizes $F(y, g_1^{[r]}, \ldots, g_{j-1}^{[r]}, u, g_{j+1}^{[r-1]}, \ldots, g_d^{[r-1]})$ over $u \in H_j$. This procedure for (6) entails the well known backfitting procedure with isotonic regressions on $X_j$, $\Pi(\cdot | H_j)$ which is given as

$$(7) \qquad g_j^{[r]} = \Pi \left( \mathbf{Y} - g_1^{[r]} - \cdots - g_{j-1}^{[r]} - g_{j+1}^{[r-1]} - \cdots - g_d^{[r-1]} \Big| H_j \right),$$

$r = 1, 2, \ldots$, $j = 1, \ldots, d$, with initial values $g_j^{[0]} = 0$ where $\mathbf{Y} = (Y^1, \ldots, Y^n)$. For a more precise description, we introduce a notation $\widetilde{Y}_j^{i,[r]} = Y^i - g_1^{i,[r]} - \cdots - g_{j-1}^{i,[r]} - g_{j+1}^{i,[r-1]} - \cdots - g_d^{i,[r-1]}$ where $g_k^{i,[r]}$ is the value of $g_k$ at $X_k^i$ after the $r$-th iteration, i.e. $\widetilde{Y}_j^{i,[r]}$ is the residual at the $j$-th cycle in the $r$-th iteration. Then, we have

$$g_j^{i,[r]} = \max_{s \leq i} \min_{t \geq i} \sum_{\ell=s}^{t} \widetilde{Y}_j^{(\ell),[r]} / (t - s + 1).$$

Here, $\widetilde{Y}_j^{(\ell),[r]}$ is the residual at the $k$-th cycle in the $r$-th iteration corresponding to the $X_j^{(\ell)}$.

Let $g^*$ be the projection of $\mathbf{Y}$ onto $H_1 + \cdots + H_d$, i.e., the minimizer of the problem (5).

**Theorem 2.** *The sequence, $g_{(r,j)} \equiv \sum_{1 \leq k \leq j} g_k^{[r]} + \sum_{j \leq k \leq d} g_k^{[r-1]}$, converges to $g^*$ as $r \to \infty$ for $j = 1, \ldots, d$. Moreover if the problem (6) has a solution that is unique* up to an additive constant, say $\mathbf{g}^* = (g_1^*, \ldots, g_d^*)$, the sequences $g_j^{[r]}$ converge to a vector with constant entries as $r \to \infty$ for $j = 1, \ldots, d$.

In general, $\mathbf{g}^*$ is not unique. Let $\mathbf{g} = (g_1, \ldots, g_d)$ be a solution of (6). Suppose, e.g. that there exists a tuple of non constant vectors $(f_1, \ldots, f_d)$ such that $f_1^i + \cdots + f_d^i = 0$ for $i = 1, \ldots, n$ and $g_j + f_j$ are monotone. Then, one does not have the unique solution for (6) since $\sum_{j=1}^d g_j^i = \sum_{j=1}^d (g_j^i + f_j^i)$ and $g_j + f_j$ are monotone. This phenomenon is similar to 'concurvity', introduced in Buja et al. (1989). One simple example for non-uniqueness is the case that elements of $X$ are ordered in the same way, i.e., $X_j^p \leq X_j^q \Leftrightarrow X_k^p \leq X_k^q$ for any $(p,q)$ and $(j,k)$. For example when $d = 2$, if $g$ solves (5), then $(\alpha g, (1-\alpha)g)$ for any $\alpha \in [0,1]$ solve (6). Other examples occur if elements of $X$ are ordered in the same way for a subregion.

## 4. Simulations

In this section, we present some simulation results for the finite sample performance. These numerical experiments are done by R on windows. We iterate 1000 times for each setting. For each iteration, we draw random samples from the following model

$$(8) \qquad Y = m_1(X_1) + m_2(X_2) + \epsilon,$$

where $(X_1, X_2)$ has truncated bivariate normal distribution and $\epsilon \sim N(0, 0.5^2)$.

In Table 1 and 2, we present the empirical MISE (mean integrated squared error) of the backfitting estimator and the oracle estimator. We also report the ratio (B/O), MISE of the backfitting estimator to MISE of the oracle estimator. For Table 1, we set $m_1(x) = x^3$ and $m_2(x) = \sin(\pi x/2)$. The results in Table 1 show that the backfitting estimator and the oracle estimator have a very similar finite sample performance. See that the ratio (B/O) is near to one in most cases and converges to one as sample size grows. We observe that when two covariates have strong negative correlation, the backfitting estimator has bigger MISE than the oracle estimator but the ratio (B/O) goes down to one as sample size grows. Figure 1 shows typical curves from the backfitting and oracle estimators for $m_1$. We show the estimators that achieve 25%, 50% and 75% quantiles of the $L_2$-distance between the backfitting and the oracle estimator for $m_1(x)$. We observe that the backfitting and the oracle estimator produce almost identical curves.

TABLE 1
*Comparison between the backfitting and the oracle estimator: Model (8) with $m_1(x) = x^3$, $m_2(x) = \sin(\pi x/2)$, sample size 200, 400, 800 and different values of $\rho$ for covariate distribution*

|     |      | $m_1$       |         |       | $m_2$       |         |       |
|-----|------|-------------|---------|-------|-------------|---------|-------|
| $n$ | $\rho$ | Backfitting | Oracle  | B/O   | Backfitting | Oracle  | B/O   |
| 200 | 0    | 0.01325     | 0.01347 | 0.984 | 0.01793     | 0.01635 | 1.096 |
|     | 0.5  | 0.01312     | 0.01350 | 0.972 | 0.01817     | 0.01674 | 1.086 |
|     | −0.5 | 0.01375     | 0.01345 | 1.022 | 0.01797     | 0.01609 | 1.117 |
|     | 0.9  | 0.01345     | 0.01275 | 1.055 | 0.01815     | 0.01601 | 1.134 |
|     | −0.9 | 0.01894     | 0.01309 | 1.447 | 0.02363     | 0.01633 | 1.447 |
| 400 | 0    | 0.00824     | 0.00839 | 0.982 | 0.01068     | 0.01000 | 1.068 |
|     | 0.5  | 0.00825     | 0.00845 | 0.977 | 0.01070     | 0.01001 | 1.063 |
|     | −0.5 | 0.00831     | 0.00830 | 1.001 | 0.01081     | 0.00997 | 1.084 |
|     | 0.9  | 0.00846     | 0.00814 | 1.040 | 0.01092     | 0.00997 | 1.095 |
|     | −0.9 | 0.10509     | 0.00805 | 1.305 | 0.01311     | 0.00992 | 1.321 |
| 800 | 0    | 0.00512     | 0.00525 | 0.976 | 0.00654     | 0.00621 | 1.053 |
|     | 0.5  | 0.00502     | 0.00513 | 0.977 | 0.00646     | 0.00614 | 1.052 |
|     | −0.5 | 0.00509     | 0.00513 | 0.994 | 0.00660     | 0.00620 | 1.066 |
|     | 0.9  | 0.00523     | 0.00500 | 1.046 | 0.00667     | 0.00611 | 1.091 |
|     | −0.9 | 0.00603     | 0.00498 | 1.211 | 0.00757     | 0.00612 | 1.220 |

Table 2 reports simulation results for the case that one component function is not smooth. Here, $m_1(x) = x, |x| > 0.5; 0.5, 0 \leq x \leq 0.5; -0.5, -0.5 \leq x < 0$ and $m_2(x) = \sin(\pi x/2)$. Even in this case the backfitting estimator shows a quite good performance. Thus the oracle property of additive isotonic least square regression

TABLE 2
*Comparison between the backfitting and the oracle estimator: Model (8) with*
$m_1(x) = x, |x| > 0.5; 0.5, 0 \leq x \leq 0.5; -0.5, -0.5 \leq x < 0$, $m_2(x) = \sin(\pi x/2)$, *sample size 200,*
*400, 800 and different values of $\rho$ for covariate distribution*

|     |      | $m_1$ |        |       | $m_2$ |        |       |
| --- | ---- | ----------- | ------ | ----- | ----------- | ------ | ----- |
| $n$ | $\rho$ | Backfitting | Oracle | B/O | Backfitting | Oracle | B/O |
| 200 | 0    | 0.01684 | 0.01548 | 1.088 | 0.01805 | 0.01635 | 1.104 |
|     | 0.5  | 0.01686 | 0.01541 | 1.094 | 0.01756 | 0.01604 | 1.095 |
|     | −0.5 | 0.01726 | 0.01541 | 1.120 | 0.01806 | 0.01609 | 1.123 |
|     | 0.9  | 0.01793 | 0.01554 | 1.154 | 0.01852 | 0.01628 | 1.138 |
|     | −0.9 | 0.02269 | 0.01554 | 1.460 | 0.02374 | 0.01633 | 1.454 |
| 400 | 0    | 0.01016 | 0.00950 | 1.071 | 0.01094 | 0.01014 | 1.079 |
|     | 0.5  | 0.00987 | 0.00944 | 1.046 | 0.01088 | 0.01025 | 1.062 |
|     | −0.5 | 0.01010 | 0.00944 | 1.070 | 0.01084 | 0.00998 | 1.086 |
|     | 0.9  | 0.01000 | 0.00897 | 1.115 | 0.01105 | 0.00996 | 1.109 |
|     | −0.9 | 0.01192 | 0.00897 | 1.330 | 0.01308 | 0.00996 | 1.314 |
| 800 | 0    | 0.00576 | 0.00552 | 1.044 | 0.00657 | 0.00622 | 1.056 |
|     | 0.5  | 0.00578 | 0.00555 | 1.041 | 0.00651 | 0.00617 | 1.055 |
|     | −0.5 | 0.00588 | 0.00555 | 1.059 | 0.00657 | 0.00614 | 1.071 |
|     | 0.9  | 0.00598 | 0.00538 | 1.110 | 0.00670 | 0.00616 | 1.088 |
|     | −0.9 | 0.00695 | 0.00538 | 1.291 | 0.00772 | 0.00612 | 1.262 |

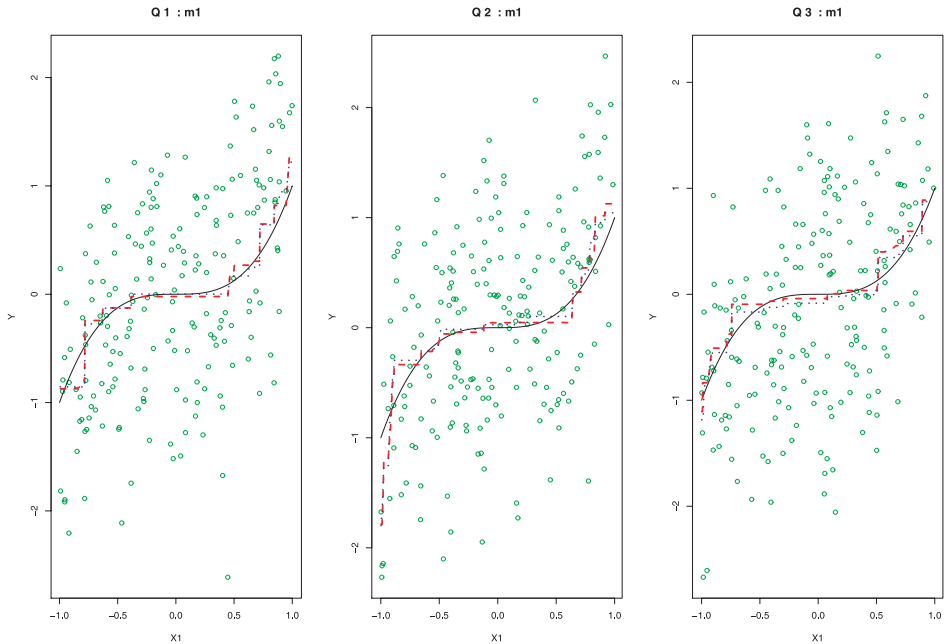

FIG 1. *The real lines, dashed lines and dotted lines show the true curve, backfitting estimates and oracle estimates, respectively. The left, center and right panels represent fitted curves for the data sets that produce 25%, 50% and 75% quantiles for the distance between the backfitting and the oracle estimator in Monte Carlo simulations with $\rho = 0.5$ and 200 observations.*

is well supported by our asymptotic theory and by the simulations.

**Appendix: Proofs**

*A.1. Proof of Theorem 1*

The proof of Theorem 1 is divided into several lemmas.

**Lemma 3.** *For $j = 1, \ldots, d$ it holds that*

$$\sup_{n^{-2/9} \leq u_j \leq 1 - n^{-2/9}} |\widehat{m}_j(u_j) - m_j(u_j)| = O_P[(\log n) n^{-2/9}].$$

*Proof.* Because $\varepsilon^i$ has subexponential tails (see (A3)) we get that $\sup_{1 \leq i \leq n} |\varepsilon^i| = O_P(\log n)$. This implies that $\max_{1 \leq j \leq d} \sup_{0 \leq u_j \leq 1} |\widehat{m}_j(u_j)| = O_P(\log n)$. We now consider the regression problem

$$Y^i/(\log n) = c/(\log n) + m_1(X_1^i)/(\log n) + \ldots + m_d(X_d^i)/(\log n) + \varepsilon^i/(\log n).$$

Now, in this model the least squares estimators of the additive components are bounded and therefore we can use the entropy bound for bounded monotone functions (see e.g. (2.6) in van de Geer [39] or Theorem 2.7.5 in van der Vaart and Wellner [40]). This gives by application of empirical process theory for least squares estimators, see Theorem 9.2 in van de Geer [39] that

$$\frac{1}{n} \sum_{i=1}^{n} \left[\widehat{m}_1(X_1^i) - m_1(X_1^i) + \ldots + \widehat{m}_d(X_d^i) - m_d(X_d^i)\right]^2 = O_P[(\log n)^2 n^{-2/3}].$$

And, using Lemma 5.15 in van de Geer [39], this rate for the empirical norm can be replaced by the $L_2$ norm:

$$\int \left[\widehat{m}_1(u_1) - m_1(u_1) + \ldots + \widehat{m}_d(u_d) - m_d(u_d)\right]^2 p(u) \, du = O_P[(\log n)^2 n^{-2/3}].$$

Because $p$ is bounded from below (see (A2)) this implies

$$\int \left[\widehat{m}_1(u_1) - m_1(u_1) + \ldots + \widehat{m}_d(u_d) - m_d(u_d)\right]^2 du = O_P[(\log n)^2 n^{-2/3}].$$

Because of our norming assumption (2) for $\widehat{m}_j$ and $m_j$ the left hand side of the last equality is equal to

$$\int \left[\widehat{m}_1(u_1) - m_1(u_1)\right]^2 du_1 + \ldots + \int \left[\widehat{m}_d(u_d) - m_d(u_d)\right]^2 du_d.$$

This gives

(9) $$\max_{1 \leq j \leq d} \int \left[\widehat{m}_j(u_j) - m_j(u_j)\right]^2 du_j = O_P[(\log n)^2 n^{-2/3}].$$

We now use the fact that for $j = 1, \ldots, d$ the derivatives $m'_j$ are bounded. This gives together with the last bound the statement of Lemma 3. □

We now define localized estimators $\widehat{m}^{OR}_{j,loc}$ and $\widehat{m}_{j,loc}$. They are defined as $\widehat{m}^{OR}_j$ and $\widehat{m}_j$ but now the sum of squares runs only over indices $i$ with $x_j - (\log n)^{1/\gamma} \times n^{-2/(9\gamma)}c_n \leq X^i_j \leq x_j + (\log n)^{1/\gamma} n^{-2/(9\gamma)} c_n$, i.e. $\widehat{m}^{OR}_{j,loc}$ minimizes

$$\sum_{i:|X^i_j - x_j| \leq (\log n)^{1/\gamma} n^{-2/(9\gamma)} c_n} \left[ m_j(X^i_j) + \varepsilon^i - \widehat{m}^{OR}_{j,loc}(X^i_j) \right]^2$$

and $\widehat{m}_{j,loc}$ minimizes

$$\sum_{i:|X^i_j - x_j| \leq (\log n)^{1/\gamma} n^{-2/(9\gamma)} c_n} \left[ Y^i - \sum_{l \neq j} \widehat{m}_l(X^i_l) - \widehat{m}_{j,loc}(X^i_j) \right]^2.$$

Here $c_n$ is a sequence with $c_n \to \infty$ slowly enough (see below). We now argue that

$$\widehat{m}_{j,loc}(x_j) = \widehat{m}_j(x_j) \text{ for } j = 1, \ldots, d \text{ and } 0 \leq x_j \leq 1$$
(10) $\qquad\qquad$ with probability tending to 1.

This follows from Lemma 3, the fact that $m_j$ fulfills (A1) and the representation (compare (4)):

$$\widehat{m}_j(x_j) = \max_{0 \leq u \leq x_j} \min_{x_j \leq v \leq 1} \frac{\sum_{i:u \leq X^i_j \leq v} \widehat{Y}^i_j}{\#\{i : u \leq X^i_j \leq v\}}, \tag{11}$$

$$\widehat{m}_{j,loc}(x_j) = \max_{x_j - (\log n)^{1/\gamma} n^{-2/(9\gamma)} c_n \leq u \leq x_j}$$
(12)
$$\min_{x_j \leq v \leq x_j + (\log n)^{1/\gamma} n^{-2/(9\gamma)} c_n} \frac{\sum_{i:u \leq X^i_j \leq v} \widehat{Y}^i_j}{\#\{i : u \leq X^i_j \leq v\}}$$

with $\widehat{Y}^i_j = Y^i - \sum_{l \neq j} \widehat{m}_l(X^i_l)$. Here $\#A$ denotes the number of elements of a set $A$.

Proceeding as in classical discussions of the case $d = 1$ one gets:

$$\widehat{m}^{OR}_{j,loc}(x_j) = \widehat{m}^{OR}_j(x_j) \text{ for } j = 1, \ldots, d \text{ and } 0 \leq x_j \leq 1$$
(13) $\qquad\qquad$ with probability tending to 1.

We now consider the functions

$$\widehat{M}_j(u_j, x_j) = n^{-1} \sum_{i:X^i_j \leq u_j} \widehat{Y}^i_j - n^{-1} \sum_{i:X^i_j \leq x_j} \widehat{Y}^i_j,$$

$$\widehat{M}^{OR}_j(u_j, x_j) = n^{-1} \sum_{i:X^i_j \leq u_j} \left[ m_j(X^i_j) + \varepsilon^i \right] - n^{-1} \sum_{i:X^i_j \leq x_j} \left[ m_j(X^i_j) + \varepsilon^i \right],$$

$$M_j(u_j, x_j) = n^{-1} \sum_{i:X^i_j \leq u_j} m_j(X^i_j) - n^{-1} \sum_{i:X^i_j \leq x_j} m_j(X^i_j).$$

For $x_j - (\log n)^{1/\gamma} n^{-2/(9\gamma)} c_n \leq u_j \leq x_j + (\log n)^{1/\gamma} n^{-2/(9\gamma)} c_n$ we consider the functions that map $\#\{i : X^i_j \leq u_j\}$ onto $\widehat{M}_j(u_j, x_j)$, $\widehat{M}^{OR}_j(u_j, x_j)$ or $M_j(u_j, x_j)$, respectively. Then we get $\widehat{m}_{j,loc}(x_j)$, $\widehat{m}^{OR}_{j,loc}(x_j)$ and $m_j(x_j)$ as the slopes of the greatest convex minorants of these functions at $u_j = x_j$.

We now show the following lemma.

**Lemma 4.** *For $\alpha > 0$ there exists a $\beta > 0$ such that uniformly for $1 \leq l, j \leq d$, $0 \leq x_j \leq 1$ and $x_j - (\log n)^{1/\gamma} n^{-2/(9\gamma)} c_n \leq u_j \leq x_j + (\log n)^{1/\gamma} n^{-2/(9\gamma)} c_n$*

$$(14) \quad \widehat{M}_j^{OR}(u_j, x_j) - M_j(u_j, x_j)$$
$$= O_P(\{|u_j - x_j| + n^{-\alpha}\}^{1/2} n^{-1/2} (\log n)^\beta),$$

$$(15) \quad \widehat{M}_j(u_j, x_j) - \widehat{M}_j^{OR}(u_j, x_j)$$
$$= -\sum_{l \neq j} n^{-1} \left[ \sum_{i: X_j^i \leq u_j} - \sum_{i: X_j^i \leq x_j} \right] \int [\widehat{m}_l(u_l) - m_l(u_l)] \, p_{X_l | X_j}(u_l | X_j^i) \, du_l$$
$$+ O_P(\{|u_j - x_j| + n^{-\alpha}\}^{2/3} n^{-13/27} (\log n)^\beta),$$

$$(16) \quad n^{-1} \left[ \sum_{i: X_j^i \leq u_j} - \sum_{i: X_j^i \leq x_j} \right] \int [\widehat{m}_l(u_l) - m_l(u_l)] \, p_{X_l | X_j}(u_l | X_j^i) \, du_l$$
$$= n^{-1} \left[ \#\{i : X_j^i \leq u_j\} - \#\{i : X_j^i \leq x_j\} \right]$$
$$\times \int [\widehat{m}_l(u_l) - m_l(u_l)] \, p_{X_l | X_j}(u_l | x_j) \, du_l$$
$$+ O_P(\{|u_j - x_j| + n^{-1}\} n^{-2\rho/(9\gamma)} (\log n)^\beta).$$

*Proof.* Claim (14) is a standard result on partial sums. Claim (16) directly follows from (A2). For a proof of claim (15) we use the following result: For a constant $C$ suppose that $\Delta$ is a difference of monotone functions on $[0, 1]$ with uniform bound $\sup_z |\Delta(z)| \leq C$ and that $Z^1, \ldots, Z^k$ is a triangular array of independent random variables with values in $[0, 1]$. Then it holds uniformly over all functions $\Delta$

$$\sum_{i=1}^k \Delta(Z^i) - E[\Delta(Z^i)] = O_P(k^{2/3}),$$

see e.g. van de Geer [39]. This result can be extended to

$$\sum_{i=1}^l \Delta(Z^i) - E[\Delta(Z^i)] = O_P(k^{2/3}),$$

uniformly for $l \leq k$ and for $\Delta$ a difference of monotone functions with uniform bound $\sup_z |\Delta(z)| \leq C$. More strongly, one can show an exponential inequality for the left hand side. This implies that up to an additional log-factor the same rate applies if such an expansion is used for a polynomially growing number of settings with different choices of $k$, $Z^i$ and $\Delta$.

We apply this result, conditionally given $X_j^1, \ldots, X_j^n$, with $Z^i = X_l^i$ and $\Delta(u) = I[n^{-2/9} \leq u \leq 1 - n^{-2/9}][\widehat{m}_l(u) - m_l(u)]/(n^{-2/9} \log n)$. The last factor is justified by the statement of Lemma 3. This will be done for different choices of $k \geq n^{1-\alpha}$. Furthermore, we apply this result with $Z^i = X_l^i$ and $\Delta(u) = \{I[0 \leq u < n^{-2/9}] + I[1 - n^{-2/9} < u \leq 1]\}[\widehat{m}_l(u) - m_l(u)]/(\log n)$ and $k \geq n^{1-\alpha}$. This implies claim (15). □

We now show that Lemma 4 implies the following lemma.

**Lemma 5.** *Uniformly for $1 \leq j \leq d$ and $n^{-1/3} \leq x_j \leq 1 - n^{-1/3}$ it holds that*

$$(17) \quad \widehat{m}_j(x_j) = \widehat{m}_j^{OR}(x_j) - \sum_{l \neq j} \int [\widehat{m}_l(u_l) - m_l(u_l)] \, p_{X_l | X_j}(u_l | x_j) \, du_l + o_P(n^{-1/3})$$

and that with a constant $c > 0$ uniformly for $1 \leq j \leq d$ and $0 \leq x_j \leq n^{-1/3}$ or $1 - n^{-1/3} \leq x_j \leq 1$

$$\widehat{m}_j(x_j) = \widehat{m}_j^{OR}(x_j) - \sum_{l \neq j} \int [\widehat{m}_l(u_l) - m_l(u_l)] \, p_{X_l | X_j}(u_l | x_j) du_l \quad (18)$$
$$+ o_P(n^{-2/9}(\log n)^c).$$

*Proof.* For a proof of (17) we use that for $n^{-1/3} \leq x_j \leq 1 - n^{-1/3}$

(19)    $\widehat{m}_j^-(x_j) \leq \widehat{m}_j(x_j) \leq \widehat{m}_j^+(x_j)$ with probability tending to 1,

(20)    $\widehat{m}_j^{OR,-}(x_j) \leq \widehat{m}_j^{OR}(x_j) \leq \widehat{m}_j^{OR,+}(x_j)$ with probability tending to 1,

(21)    $\sup_{0 \leq x_j \leq 1} \widehat{m}_j^{OR,+}(x_j) - \widehat{m}_j^{OR,-}(x_j) = o_P(n^{-1/3})$,

where

$$\widehat{m}_j^-(x_j) = \max_{x_j - e_n \leq u \leq x_j - d_n} \min_{x_j \leq v \leq x_j + e_n} \frac{\sum_{i: u \leq X_j^i \leq v} \widehat{Y}_j^i}{\#\{i : u \leq X_j^i \leq v\}},$$

$$\widehat{m}_j^+(x_j) = \max_{x_j - e_n \leq u \leq x_j} \min_{x_j + d_n \leq v \leq x_j + e_n} \frac{\sum_{i: u \leq X_j^i \leq v} \widehat{Y}_j^i}{\#\{i : u \leq X_j^i \leq v\}},$$

$$\widehat{m}_j^{OR,-}(x_j) = \max_{x_j - e_n \leq u \leq x_j - d_n} \min_{x_j \leq v \leq x_j + e_n} \frac{\sum_{i: u \leq X_j^i \leq v} m_j(X_j^i) + \varepsilon^i}{\#\{i : u \leq X_j^i \leq v\}},$$

$$\widehat{m}_j^{OR,+}(x_j) = \max_{x_j - e_n \leq u \leq x_j} \min_{x_j + d_n \leq v \leq x_j + e_n} \frac{\sum_{i: u \leq X_j^i \leq v} m_j(X_j^i) + \varepsilon^i}{\#\{i : u \leq X_j^i \leq v\}},$$

compare (11) and (12). Here, $e_n = (\log n)^{1/\gamma} n^{-2/(9\gamma)} c_n$ and $d_n$ is chosen as $d_n = n^{-\delta}$ with $1/3 < \delta < 4/9$. Claims (19) and (20) follow immediately from the definitions of the considered quantities and (10) and (13). Claim (21) can be established by using well known properties of the isotone least squares estimator. Now, (15),(16),(19) and (20) imply that uniformly for $1 \leq j \leq d$ and $n^{-1/3} \leq x_j \leq 1 - n^{-1/3}$

$$\widehat{m}_j^{\pm}(x_j) = \widehat{m}_j^{OR,\pm}(x_j) - \sum_{l \neq j} \int [\widehat{m}_l(u_l) - m_l(u_l)] \, p_{X_l | X_j}(u_l | x_j) du_l + o_P(n^{-1/3}).$$

This shows claim (17).

For the proof of (18) one checks this claim separately for $n^{-7/9}(\log n)^{-1} \leq x_j \leq n^{-1/3}$ or $1 - n^{-1/3} \leq x_j \leq 1 - n^{-7/9}(\log n)^{-1}$ (case 1) and for $0 \leq x_j \leq n^{-7/9}(\log n)^{-1}$ or $1 - n^{-7/9}(\log n)^{-1} \leq x_j \leq 1$ (case 2). The proof for Case 1 is similar to the proof of (17). For the proof in Case 2 one considers the set $I_{j,n} = \{i : 0 \leq X_j^i \leq n^{-7/9}(\log n)^{-1}\}$. It can be easily checked that with probability tending to 1 it holds that $n^{-2/9} \leq X_l^i \leq 1 - n^{-2/9}$. Therefore it holds that $\sup_{i \in I_{j,n}} |\widehat{m}_l(X_l^i) - m_l(X_l^i)| = O_P[(\log n) n^{-2/9}]$, see Lemma 3. Therefore for $0 \leq x_j \leq n^{-7/9}(\log n)^{-1}$ the estimators $\widehat{m}_j(x_j)$ and $\widehat{m}_j^{OR}(x_j)$ are local averages of quantities that differ by terms of order $O_P[(\log n) n^{-2/9}]$. Thus also the difference of $\widehat{m}_j(x_j)$ and $\widehat{m}_j^{OR}(x_j)$ is of order $O_P[(\log n) n^{-2/9}]$. This shows (18) for Case 2.   □

We now show that Lemma 5 implies the statement of the theorem.

*Proof of Theorem 1.* We rewrite equations (17) and (18) as

$$\widehat{m} = \widehat{m}^{OR} + H(\widehat{m} - m) + \Delta, \tag{22}$$

where $\widehat{m}$, $\widehat{m}^{OR}$ and $\Delta$ are tuples of functions $\widehat{m}_j$, $\widehat{m}_j^{OR}$ or $\Delta_j$, respectively. For $\Delta_j$ it holds that

$$\sup_{n^{-1/3} \leq x_j \leq 1 - n^{-1/3}} |\Delta_j(x_j)| = o_P(n^{-1/3}), \tag{23}$$

$$\sup_{0 \leq x_j \leq 1} |\Delta_j(x_j)| = o_P(n^{-2/9}(\log n)^c). \tag{24}$$

Furthermore, $H$ is the linear integral operator that corresponds to the linear map in (17) and (18). For function tuples $f$ we denote by $Nf$ the normalized function tuple with $(Nf)_j(x_j) = f_j(x_j) - \int f_j(u_j)du_j$ and we introduce the pseudo norms

$$\|f\|_2^2 = \int [f_1(x_1) + \ldots + f_d(x_d)]^2 p(x)\, dx,$$

$$\|f\|_\infty = \max_{1 \leq j \leq d} \sup_{0 \leq x_j \leq 1} |f_j(x_j)|.$$

Here $p_j$ is the marginal density of $X_j^i$ and $p$ is the joint density of $X^i$. We make use of the following properties of $H$. On the subspace $\mathcal{F}_0 = \{f : f = Nf\}$ the operator $H$ has bounded operator norm:

$$\sup_{f \in \mathcal{F}_0, \|f\|_2 = 1} \|Hf\|_2 = O(1). \tag{25}$$

For the maximal eigenvalue $\lambda_{\max}$ of $H$, it holds that

$$\lambda_{\max} < 1. \tag{26}$$

Claim (25) follows directly from the boundedness of $p$. Claim (26) can be seen as follows. Compare also with Yu, Park and Mammen [42].

A simple calculation gives

$$\int (m_1(u_1) + \cdots + m_d(u_d))^2 p(u)\, du = \|m\|_2^2 = \int m^T(I - H)m(u)p(u)du. \tag{27}$$

Let $\lambda$ be an eigenvalue of $H$ and $m_\lambda$ be an eigen(function)vector corresponding to $\lambda$. With (27), we have

$$\|m_\lambda\|_2^2 = \int m_\lambda^T(I - H)m_\lambda(u)p(u)du = (1 - \lambda)\int m_\lambda^T m_\lambda(u)p(u)du.$$

Thus, the factor $1 - \lambda$ must be strictly positive, i.e. $\lambda < 1$. This implies $I - H$ is invertible and hence we get that

$$N(\widehat{m} - m) = (I - H)^{-1}N(\widehat{m}^{OR} - m) + (I - H)^{-1}N\Delta.$$

Here we used that because of (22)

$$N(\widehat{m} - m) = N(\widehat{m}^{OR} - m) - NH(\widehat{m} - m) + N\Delta$$
$$= N(\widehat{m}^{OR} - m) - HN(\widehat{m} - m) + N\Delta.$$

We now use
$$(I - H)^{-1} = I + H + H(I - H)^{-1}H,$$
$$(I - H)^{-1} = I + H(I - H)^{-1}.$$

This gives
$$N(\widehat{m} - m) = N(\widehat{m}^{OR} - m) + N\Delta + HN(\widehat{m}^{OR} - m)$$
$$+ H(I - H)^{-1}HN(\widehat{m}^{OR} - m) + H(I - H)^{-1}\Delta.$$

We now use that

(28) $$\|HN(\widehat{m}^{OR} - m)\|_2 \leq \|HN(\widehat{m}^{OR} - m)\|_\infty = o_P(n^{-1/3}),$$

(29) $$\sup_{f \in \mathcal{F}_0, \|f\|_\infty = 1} \|Hf\|_\infty = O(1).$$

Claim (28) follows because $\widehat{m}^{OR}$ is a local average of the data, compare also Groeneboom [12], Groeneboom, Lopuhaa and Hooghiemstra [14] and Durot [8]. Claim (29) follows by a simple application of the Cauchy Schwarz inequality, compare also (85) in Mammen, Linton and Nielsen [28].

This implies that
$$\|N(\widehat{m} - m) - N(\widehat{m}^{OR} - m) - N\Delta\|_\infty = o_P(n^{-1/3}).$$

Thus,
$$\sup_{n^{-1/3} \leq x_j \leq 1 - n^{-1/3}} |N(\widehat{m} - \widehat{m}^{OR})_j(x_j)| = o_P(n^{-1/3}),$$
$$\sup_{0 \leq x_j \leq 1} |N(\widehat{m} - \widehat{m}^{OR})_j(x_j)| = o_P(n^{-2/9}(\log n)^c)$$

This implies the statement of Theorem 1. □

### A.2. Proof of Theorem 2

For a given closed convex cone $K$, we call $K^* \equiv \{f : \langle f, g \rangle \leq 0 \text{ for all } g \in K\}$ the dual cone of $K$. It is clear that $K^*$ is also a convex cone and $K^{**} = K$. It is pointed out in Barlow and Brunk [3] that if $P$ is a projection onto $K$ then $I - P$ is a projection onto $K^*$ where $I$ is the identity operator. Let $P_j$ be a projection onto $H_j$ then $P_j^* \equiv I - P_j$ is a projection onto $H_j^*$. The backfitting procedure (7) to solve the minimization problem (5) corresponds in the dual problem to an algorithm introduced in Dykstra [9]. See also Gaffke and Mathar [11]. We now explain this relation. Let $H_j$, $j = 1, \ldots, d$, be sets of monotone vectors in $\mathbb{R}^n$ with respect to the orders of $X_j$ and $P_j = \Pi(\cdot|H_j)$. Denote the residuals in algorithm (7) after the $k$-th cycle in the $r$-th iteration with $h_{(r,k)}$. Then, we have

$$h_{(1,1)} = \mathbf{Y} - g_1^{[1]} = P_1^*\mathbf{Y},$$
$$h_{(1,2)} = \mathbf{Y} - g_1^{[1]} - g_2^{[1]} = P_1^*\mathbf{Y} - P_2 P_1^*\mathbf{Y} = P_2^* P_1^*\mathbf{Y},$$
$$\vdots$$
(30) $$h_{(1,d)} = \mathbf{Y} - g_1^{[1]} - \cdots - g_d^{[1]} = P_d^* \cdots P_1^*\mathbf{Y};$$

$$h_{(r,1)} = \mathbf{Y} - g_1^{[r]} - g_2^{[r-1]} - \cdots - g_d^{[r-1]} = P_1^*(\mathbf{Y} - g_2^{[r-1]} - \cdots - g_d^{[r-1]}),$$

$$\vdots$$

$$h_{(r,k)} = \mathbf{Y} - g_1^{[r]} - \cdots - g_k^{[r]} - g_{k+1}^{[r-1]} - \cdots - g_d^{[r-1]}$$
$$= P_k^*(\mathbf{Y} - g_1^{[r]} - \cdots - g_{k-1}^{[r]} - g_{k+1}^{[r-1]} - \cdots - g_d^{[r-1]}),$$

$$\vdots$$

(31) $\quad h_{(r,d)} = \mathbf{Y} - g_1^{[r]} - \cdots - g_d^{[r]} = P_d^*(\mathbf{Y} - g_1^{[r]} - \cdots - g_{d-1}^{[r]}).$

With the notation $I_{r,k} \equiv -g_k^{[r]}$ for the incremental changes at the $k$-th cycle in the $r$-th iteration, equations (30) and (31) form a Dykstra algorithm to solve the following optimization problem:

(32) $$\min_{h \in H_1^* \cap \cdots \cap H_d^*} \sum_{i=1}^n (Y^i - h^i)^2.$$

Denote the solutions of (32) with $h^*$. Theorem 3.1 of Dykstra [9] shows that $h_{(r,j)}$ converges to $h^*$ as $r \to \infty$ for $j = 1, \ldots, d$. From the dual property, it is well known $g^* = \mathbf{Y} - h^*$ and also it is clear that $g_{(r,j)} = \mathbf{Y} - h_{(r,j)}$ for $j = 1, \ldots, d$. Since $h_{(r,j)}$ converges to $h^*$, $g_{(r,j)}$ converge to $g^*$ as $r \to \infty$ for $j = 1, \ldots, d$. The convergence of $g_j^{[r]}$ follows from Lemma 4.9 of Han [17].